%% file: main.tex
\documentclass{sig-alternate}

\newtheorem{theorem}{Theorem}
\newtheorem{definition}{Definition}
\newtheorem{proposition}{Proposition}
\newtheorem{observation}{Observation}

\newcommand{\cqfd}{\mbox{}\nolinebreak\hfill\rule{2mm}{2mm}\medbreak\par}
\usepackage{algpseudocode}
\usepackage[linesnumbered,ruled,vlined]{algorithm2e}
\usepackage{caption}
\usepackage{subcaption}

\title{Revenue maximization with access and information pricing schemes in a partially-observable queueing game}

\numberofauthors{2}
\author{Tesnim Naceur$^\dagger$ and Yezekael Hayel$^\dagger$\\{~}
{\small $^\dagger$LIA/CERI, University of Avignon, France} \\{~}
{\small E-mail: \{tesnim.naceur,yezekael.hayel\}@univ-avignon.fr}\\{~}
}

\begin{document}
\maketitle
\begin{abstract}
Today's queueing network systems are more rapidly evolving and more complex than those of even a few years ago. The goal of this paper is to study customers' behavior in an unobservable Markovian M/M/1 queue where consumers have to choose between two strategic decisions about information acquisition before joining or not the queue. According to their decision, customers decide to give up the service (balk the system) or to join the queue. We study the Nash equilibrium strategies: we compute the equilibrium and we prove its existence and uniqueness. Based on this result, we consider the problem of revenue maximization where the provider has to choose between two charging mechanisms: to charge the access to the system or to charge the queue length information to new incoming customers. We propose an heuristic algorithm to solve the considered problem and numerical experiments have been conducted in order to illustrate the result. We show that, depending on the sensitivity of customers to their waiting time, the provider will change the pricing policy used. Particularly, when customers are more sensitive then the provider will charge the information instead of the access.  \\
\end{abstract}

\vspace{0.5cm}
\section{Introduction}
\input{introduction}

\vspace{0.5cm}
\section{Model description}
\input{model}

\vspace{0.5cm}
\section{Revenue maximization}
\input{revenue}

\vspace{0.5cm}
\section{Numerical Illustrations}
\input{numerical}

\vspace{0.5cm}
\section{Conclusions}
In this work we have studied an unobservable queueing game in which customers decide to get or not the information about actual queue length in order to make their joining decision. This strategic decision depends on numerous parameters of the system. In a first step, we have been able to prove that this queueing game has a unique equilibrium and we have been able to give his closed-form expression. In a second step, we have analyzed the optimal pricing decision for the provider in order to maximize his revenue, considering that the provider has to choose between charging the information or the access to new customers. Interestingly, we have obtained that charging the access of the service is not always the best pricing policy, and this result depends on the waiting cost. Particularly, when customers are very impatient and strategic, it is preferable for the provider to charge the information instead of the access to the service.

\bibliography{mybib}

\appendix
\input{appendix}

\end{document}

%% file: introduction.tex
Waiting in unobservable queue systems causes frustration and annoyance to customers. The good news is that customers can make the situation better by acquiring on-line information of the queue-length and decide after that, whether to join the system or balk. The cost to acquire this information can be very different depending on the context. In many cases, this information is available at no charge, but still it is not costless to the customers who spend time and effort to obtain it, whether they pre-register, download relevant application, or even acquire relevant terminal interface equipment. 

There are numerous examples of queuing systems which provide the queue-length information on-line. This is the case, for example, of Google restaurant requests, where the average wait times are now viewable in Maps to approximately one million restaurants around the world. Devices follow the waiting time and inform searchers how long they should expect to wait. Also, many hospitals and health-care services in the United States and Canada have recently provided current waiting times on its web site. For Online queue-length information in the hospitals you can see: Florida Hospital \footnote{www.floridahospital.com}, The Ontario site web Reston hospital center \footnote{www.restonhospital.com} and the JFK medical center \footnote{jfkmc.com/our-services/er-wait-time.dot}. 
In our problem, we consider customers arriving to a single-server queue. The decision taken by each customer is not related to joining or not the system, but to acquire or not the information. Customers may choose to acquire information about the queue length with a fixed cost in order to take their decision. The cost of inspection is imposed by the service provider. If a customer inspects the queue, the next step is similar to Naor \cite{Naor69} observable-queue model. Otherwise, customer decides to join the system without observation.\\
In our model, we prove the existence and uniqueness of the equilibrium and we analyze the impact of inspection cost and the access fee on the revenue of the service provider.
The study of queueing systems was started with Naor  \cite{Naor69} who studied the observable M/M/1 queue and assumed that an arriving customer observes the number of customers and then makes his decision whether to join or balk the queueing system. Then, Edelson and Hildebrand \cite{Edelson} considered the same problem in unobservable case: they assumed that the customers make their decisions without being informed about the queue-length information. Articles \cite{Shone13,simhon2016optimal} was compared between extreme configurations when all costumers have queue-length information or no one have it. Recently, R. Hassin and R. Roet-Green \cite{Hassin17} have showed a deep and comprehensive analysis for the queueing M/M/1 model where the costumers choose to: join the queue, balk, or inspect the queue length. They have computed the equilibrium and analyzed the impact of inspection cost on equilibrium, revenue and social welfare. They discussed both homogeneous and heterogeneous customers, and obtained various results. However, the authors do not have a closed-form expression of the equilibrium strategy. Our model can be seen as a special case of theirs (restricting the decision space to 2 actions) but, by doing that, we are able to study a high level optimization problem for the provider about which pricing mechanism to use in order to maximize his revenue. Information structure proposed to customers is studied in \cite{Allon11}. Particularly, the authors study how the service provider can use the delay announcements, like estimates of waiting time or place in the queue, in order to induce the appropriate customers' behavior, so as to maximize her own profits. Recent empirical studies have been performed in hospital networks in United States. Those studies have shown that oscillations in the incoming arrival of patients can be observed when patients are sensitive to delay announcements \cite{Dong18}. These oscillations may cause an increase in the average waiting times for patients. This paper provides empirical evidence that patients take delay information into account for making their decision. A recent survey \cite{Hassin16} describes almost all the research articles related to strategic queueing systems. In this survey book, chapter 3 is related to information in these problems. And, as a reader can see, there are no articles that study the problem in which the strategic decision of each customer is to acquire or not the information; and not to join or balk the queue. Many articles study the impact of heterogeneity in the population between informed and uniformed customers, design of optimal signaling mechanisms and processing time information. In contrast with previous papers, our work focus on the information acquisition decision process by each customer, and how this decision is affected by provider's cost and how it affects the decision of every customers (equilibrium situation). The main contributions of our work are:
\begin{itemize}
\item to model an information acquisition game for an M/M/1 queue with strategic customers,
\item to characterize the equilibrium situation in which no customers has an interest to deviate about his decision to acquire or not the information,
\item to study the impact of information cost on a global performance which is the revenue of the service provider,
\item to compare two charging mechanisms in order to optimize the provider's revenue: to charge the access to the queue or the information to customers. 
\end{itemize}

This paper is organized as follows: section 2 describes the mathematical model and the equilibrium analysis. Section 3 deals with the revenue maximization problem and in section 4 we analyze analyze the best pricing mechanism between charging the access to the queue or the information to customers. Our results are numerically illustrated in section 5 and section 6 concludes our work.

All proofs are given in the appendix.

%% file: model.tex
\begin{figure}[h]
\centering
\includegraphics[scale=0.22]{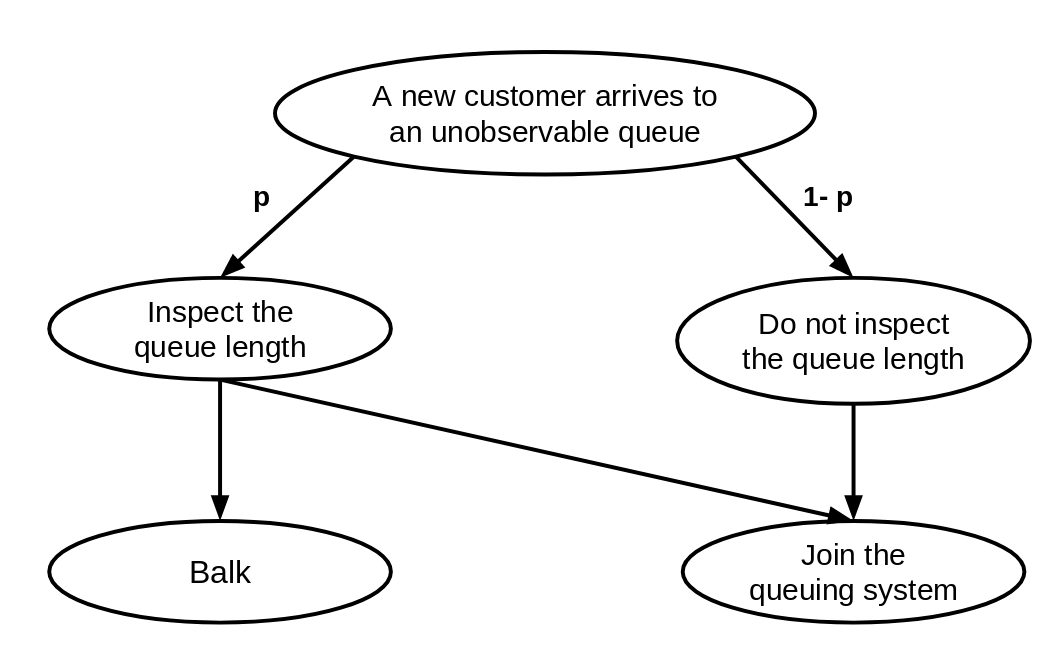} 
\caption{\label{fig:decision1}The decision process of each arriving customer. $p$ denotes the probability that a new customer acquires the queue-length information.}
\end{figure}

We consider an unobservable $M/M/1$ queue system where customers arrive according to a Poisson process with intensity $\lambda$. Service times are assumed to be independent and exponentially distributed with rate $\mu$. We assume arriving customers face two possible actions : to acquire the queue-length information or not. Based on the action taken, each customer decides to join or balk the system. Figure \ref{fig:decision1} shows the customer's decision process. The information acquisition process is out of scope of the work and we assume that every customer gets perfect information in real time, meaning that it is exactly the actual number of customers waiting in the system. The cost of each time unit spent in the system is noted by $C_W$ and the cost of inspection, imposed by the provider, is $C_I$ where $C_W >$ 0 and $C_I \ge $ 0. If it is not directly charged by the service provider, the latter inspection cost may incur from the information acquisition process that needs specific tools (an access to the Internet, to call a data center, to reach a specific location, etc.). Finally, all customers who join the system have the same service valuation $R>0$. Note that, like in most papers related to strategic decisions in queueing systems, we consider homogeneous population of customers , i.e. same service time and costs. The generalization to heterogeneous customers is proposed as a perspective of this work. 
\vspace{0.5cm}

\subsection{Joining process}
The customer who gets the information, decides to join the queue or balk based on queue-length information and evaluation of the service. According to Naor's model, a customer in an observable queue system with $i$ customers already in the queue, joins if and only if $i \leq n_e - 1$, where: 
\begin{eqnarray}
n_e =\lfloor \frac{R\mu}{C_W } \rfloor.
\label{Naor}
\end{eqnarray}
The system utilization factor is defined by: $\rho= \frac{\lambda}{\mu}$, where $\lambda < \mu $. A consequence of this joining process is that the effective arrival rate $\lambda_e$ of customers to the system depends on the current number of customers in the queue, that is:
\begin{equation*}
\lambda _e   = \left\{
    \begin{array}{ll}
        \lambda & \mbox{if } i \leq n_e-1,  \\
        (1-p)\lambda & \mbox{otherwise.}
    \end{array}
\right.
\end{equation*}
We look now for determining the stationary distribution of the Markov process considering state depending rates. The balance equation for states $i \in \{ 0, \ldots, n_e -1\}$ is given by:
\begin{equation*}
\lambda \pi _i = \mu \pi _{i+1},
\end{equation*}
where $\pi_i$ is the stationary probability in state $i$. For states $ i \geq n_e$, the balance equation is:
\begin{equation*}
(1-p)\lambda \pi_i= \mu \pi _{i+1}.
\end{equation*}
To simplify the following computations, we define $\eta = 1-(1-p)\rho$. Then, the stationary distribution for any states $i > 0 $, depending on stationary probability $\pi_0$, is:
\begin{equation*}
\begin{split}
\pi_i   = \left\{
    \begin{array}{ll}
        \rho^i \pi_0 & \mbox{if } i \leq n_e-1,  \\
        \rho^{n_e}(1-\eta)^{i-n_e}\pi_0 & \mbox{otherwise.}\\
    \end{array}
\right.
\end{split}
\end{equation*}
Finally, by summing to one, the stationary probability $\pi_0$ that the queueing system has no customer is: 
\begin{equation*}
\begin{split}
\pi_0 &=[ \sum_{i=0}^{n_e-1} \rho^i +\sum_{i=n_e}^{\infty} (1-\eta)^{i-n_e}]^{-1}=[\frac{ 1- \rho^{n_e}}{1-\rho}+ \frac{ \rho^{n_e} }{\eta }]^{-1}.  
\end{split}
\end{equation*}
The stationary probability $\pi_0$ that the queueing system has no customer is: 
\begin{equation*}
\pi_0= \frac{ \eta }{n_e\eta +1}.
\end{equation*}
The expected utility $U_I$ for a customer who gets the queue length information is given by:
\begin{eqnarray}
U_I&=& \sum_{i=0}^{n_e-1}  \pi_i(R-C_W\frac{ i+1 }{\mu })- C_I.
\end{eqnarray}
The expected utility for getting the queue length information can be simplified to:
\begin{equation*}
\begin{split}
U_I   &=\pi_0[R\frac{ 1 -\rho^{n_e} }{1-\rho }-\frac{ C_W }{\mu }  \frac{ 1-(n_e+1)\rho^{n_e}+n_e\rho^{n_e+1} }{(1-\rho)^2 }]-C_I.
\end{split}
\end{equation*}
The expected utility $U_{NI}$ for a customer who do not get the queue length information and joins directly the system is given by: 
\begin{eqnarray}
U_{NI}&=&\sum_{i=0}^{\infty} \pi_i (R-C_W\frac{ i+1 }{\mu }).
\end{eqnarray}
This expected utility for not getting the queue length information can be simplified to the following expression:
\begin{eqnarray*}
U_{NI}&=&(R-\frac{ C_W }{\mu }\pi_0 [\frac{ 1-(n_e+1)\rho^{n_e}+n_e\rho^{n_e+1} }{(1-\rho)^2 } \\
&&+ \frac{ \rho^{n_e}(n_e\eta+1) }{\eta^2 } ]).
\end{eqnarray*}

We now look for equilibrium in terms of proportion of customers who decide to get the queue length information. Note that this proportion is equal to the probability with which each customer will decide to get the queue length information. Indeed, each customer takes a deterministic action (to get the queue length information or not) with the same probability distribution, that is here a Bernoulli distribution with parameter $p$. Therefore, the proportion of customers who get the queue length information is exactly $p$. Before investigating the equilibrium, let us give the first property of the expected utilities depending on this proportion $p$. 
\begin{observation} \label{ob1}
We consider that there is no cost to get the queue length information, i.e. $C_I$ = 0, then:
$$
\forall p \in [0,1], \quad U_{NI}(p)< U_I(p).
$$
\end{observation}

This observation states an intuitive result that implies that any customer will decide to get the queue length information when this action does not imply any specific additional cost. A corollary of this observation is that, when $C_I=0$, the equilibrium is $p^*=1$.

\subsection{Equilibrium}

Previously, we have been able to demonstrate a first simple intuitive result on the equilibrium in a very specific context where there is no cost for information. We look now for an equilibrium in the general case. Let us first define an equilibrium for our queueing game with indistinguishable infinitely many players. In this case, we can study the equilibrium by considering the utility of one tagged customer against the whole population. Therefore, we denote by $U(p,p')$ the (expected) utility of customer who decides to get the queue length information with probability $p$ when everyone else get this information with probability $p'$. Then, we have:
$$
U(p,p')=pU_{I}(p')+(1-p)U_{NI}(p').
$$
We can now formalize an equilibrium strategy as follows.
\begin{definition}
A (symmetric) equilibrium is a strategy $p^*$ such that:
$$
p^* =\arg\max_{p}U(p,p^*).
$$
\end{definition}
An equilibrium is a best response against itself and equilibrium here is in the sense of Nash equilibrium. Based on the equilibrium formal definition, we are able in the following proposition to express explicitly the proportion of customers who decide to get the queue length information depending on the parameter of the system.\\
To simplify the presentation, we define the following variables:\\
$$
K_1=\rho^{n_e}( C_W(n_e+1)-R\mu),
$$
$$
K_2=C_I(1-\rho^{n_e+1}),
$$
$$
K_3=\mu C_I(1-\rho^{n_e})\rho^2,
$$
\begingroup\makeatletter\def\f@size{8}\check@mathfonts
$$
K_4=-\left(C_I(1-\rho^{n_e})\rho \mu+K_2 \rho \mu+(1-\rho)(\rho^{n_e+1}C_W-K_1 \rho) \right),
$$
\endgroup
$$
K_5=K_2 \mu-K_1 (1-\rho),
$$
$$
\Delta= K_4^2-4K_3K_5.
$$

\begin{theorem}\label{theo1}
For each parameters set of the system, an equilibrium $p^*$ exists and is uniquely defined as:
\begingroup\makeatletter\def\f@size{8}\check@mathfonts
$$
p^*=\frac{2\rho-1}{2\rho}-\frac{(1-\rho^{n_e+1})}{2\rho(1-\rho^{n_e})}-\frac{(1-\rho)(\rho^{n_e+1}C_W-K_1\rho)-\sqrt \Delta}{2\mu C_I(1-\rho^{n_e})\rho^2},
$$
if $p^* \in [0,1]$. Else, the equilibrium is either $0$ or $1$.
\endgroup
\end{theorem}

%% file: revenue.tex
In this section, we investigate the problem of maximizing the provider's revenue. This type of problem is mostly answered by considering charging the access to the queueing system. The most well-known problem is the optimal toll pricing mechanism which is applied mainly in routing games \cite{tollpricing}. In our problem, the information given to customers affect their behavior and therefore a provider may optimize his revenue by charging the information instead of the access to the queue. And, we can expect that the maximum of the revenue obtained by the provider can be higher by using such charging mechanism compared to the one obtained by charging the access as a toll. The rest of the paper aims to answer this question: is it more profitable for the provider to charge the information instead of the access to the system? And we will see that the answer is highly dependent of the customers sensitivity to their waiting time.

The objective function, provider's revenue, is the same considering both charging mechanism but the provider's control is different. We denote by $R_I$ (resp. $R_A$) the provider's revenue that depends on the information price $C_I$ (resp. the access price $C_{acc}$). Note that the variable $C_I$ that was denoted as a cost in previous section is now controlled and becomes a price/charge imposed by the provider. Then, the objective function problem for the provider considering the information charging mechanism can be written as:
\begin{eqnarray}
R_I(C_I)=\lambda p^*(C_I) C_I,
\end{eqnarray}
where $p^*(C_I)$ is the proportion of customers who accept to pay for the information, and
\begin{eqnarray}
R_A(C_{Acc})=\lambda q^*(C_{Acc}) C_{Acc},
\end{eqnarray}
where $q^*(C_{Acc})$ is the proportion of customers who pay the access to the queue. Note that in this context, customers who do not pay, do not have access to the service whereas, in the other case, where the service provider charge the queue length information, some customers do not join even if they have paid for getting the information.
\vspace{0.5cm}
\subsection{Charging the access to the queue}

If the provider decides to charge the access to the queue, its optimization problem can be written as:
\begin{eqnarray*}
\max_{C_{Acc}}R_A(C_{Acc})=\max_{C_{Acc}}\lambda q^*(C_{Acc}) C_{Acc},
\end{eqnarray*}
considering that depending on the access price $C_{Acc}$, part of the customers join the queue and part of then balk; given that customers do not know the queue length. This problem corresponds to the unobservable strategic queueing problem studied in \cite{Edelson}.  The equilibrium (proportion of customers who join the queue) depending on the parameters of the system and particularly the access fee is given by:\\

\begin{equation*}
q^*(C_{Acc}) =
 \left\{
    \begin{array}{ll}
        1 & \mbox{ if    } C_{Acc} \le R-\frac{C_W}{\mu-\lambda},  \\
        \frac{\mu-\frac{C_W}{R-C_{Acc}}}{\lambda} & \mbox{ if    } R-\frac{C_W}{\mu-\lambda}  \le C_{Acc} \le R-\frac{C_W}{\mu}, \\
        0 & \mbox{ if   }  C_{Acc} > R-\frac{C_W}{\mu}.
    \end{array}
\right.
\end{equation*}
We observe that the equilibrium $q^*(C_{Acc})$ is decreasing with the access fee $C_{Acc}$ which is an intuitive result. Moreover, as we can expect, the revenue is zero if the access fee is zero and, at the opposite, when the access fee goes to infinity the revenue tends to be also zero. Then, by continuity of the function, there exists almost one optimum access fee $C_{Acc}^*$ that maximizes the provider revenue $R_A(C_{Acc})$ when the provider decides to charge the queue access to the customers. In fact, we can prove easily that this maximum is unique and we denote it by $R_A^*:=\max_{C_{Acc}}R_A(C_{Acc})$.

\begin{proposition}\label{prop2}
When the provider decides to charge the queue access to customers, the maximum revenue $R_A^*$ is given by:
$$
R_A^*=\max\{R_A(R-\frac{C_W}{\mu-\lambda}), R_A(R-\sqrt{\frac{RC_W}{\mu}})\}.
$$
\end{proposition}

Then, the optimal access fee is either $R-\frac{C_W}{\mu-\lambda}$ or $R-\sqrt{\frac{RC_W}{\mu}}$.
In next section, we study the optimization of the provider's revenue when the provider decides to charge the information instead of the queue access to customers. 
\vspace{0.5cm}

\subsection{Charging the information}
Considering that the provider is charging the information and not the access, it means that the control variable of the provider is the parameter $C_I$. Then, the provider determines the optimal information fee $C_I^*$ that optimize his revenue which is expressed by:
$$
R_I(C_I)=\lambda p^*(C_I) C_I,
$$
where $p^*(C_I)$ is the equilibrium (proportion of customers who pay for the information queue length) that depends on the information fee.

\begin{figure}[h]
\centering
\includegraphics[scale=0.25]{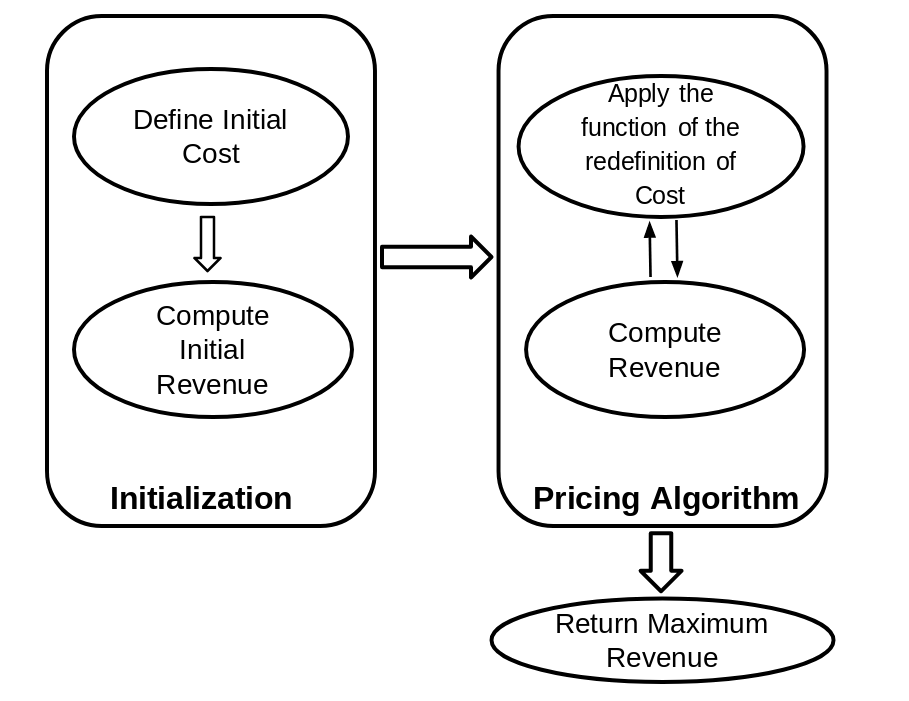} 
\caption{\label{fig:algo} Flowchart of Revenue Optimization Heuristic.}
\end{figure}

We propose an heuristic based on the redefinition of the inspection cost/access fee at each iteration, which influences the equilibrium and then the value of the revenue. The algorithm follows the process described in figure \ref{fig:algo}.

In fact, due to the strict concavity property of the revenue, the algorithm starts the exploration with a price equals to zero and increases with a fixed step this value until the revenue decreases. At that final step, the algorithm has reached the maximum value of the revenue denoted by $R_I^*$.

Based on the optimization results of the two previous pricing mechanisms, the provider is able to determine the best mechanism depending on an important parameter which is the waiting cost $C_W$. Indeed, this parameter describes how customers are sensitive to waiting time duration into the queue. 

\vspace{0.5cm}
\section{Optimal pricing policy}

\begin{figure}[h]
\centering
\includegraphics[scale=0.27]{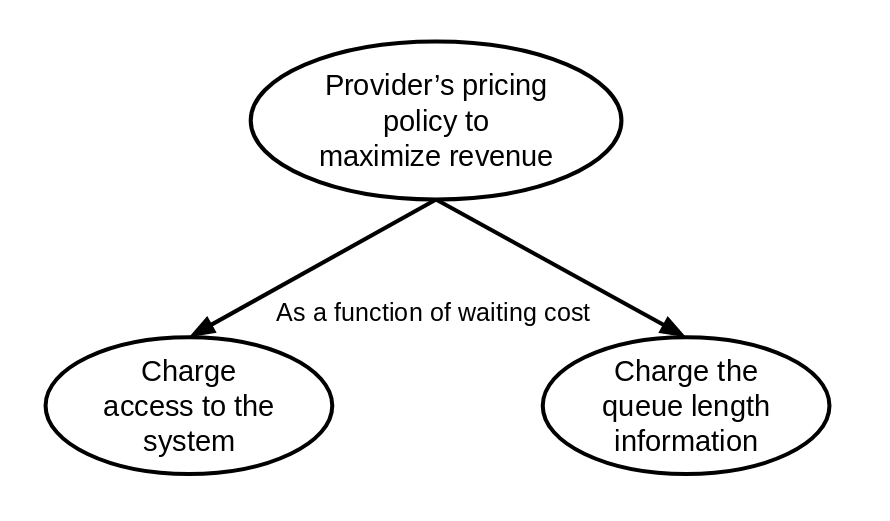} 
\caption{\label{fig:price}Choosing the best pricing mechanism for the provider.}
\end{figure}

We assume that the provider wants to determine which pricing mechanism between charging access or information leads to the best revenue. Particularly, we are interested in the sensitivity of this choice depending on an important social parameter which is the waiting cost $C_W$ of customers into the queue. Such problem is depicted on figure (\ref{fig:price}).

First we demonstrate that optimal revenue is strictly decreasing as a function of waiting cost when the provider charge the access to the queue and when he charges the queue length information, optimal revenue is equal to zero when $C_W$ tends to zero and goes to infinity when $C_W$ tends to infinity. And second, as a corollary of these results, we show that when the waiting cost is low (i.e. customers are not very sensitive to the waiting time) the provider should charge the access, whereas when the waiting cost is high the provider should preferably charge the information. Note that the estimation of the waiting cost can be obtained through statistical analysis of customers experience data extracted for surveys. This problem is related to numerous works that have been done considering impatience in customers decisions in such complex strategic queueing systems \cite{queuebook}.

\begin{proposition}\label{prop3}
When the service provider charges the access to the queue, its optimum revenue is strictly decreasing with the waiting cost.
\end{proposition}

We look now to the behavior of the optimum revenue when the service provider decides to charge the information given to incoming customers.
\begin{proposition} \label{prop4}
When the service provider charges the information, its optimum revenue is equal to zero when $C_W$ tends to zero and goes to infinity when $C_W$ tends to infinity.
\end{proposition}

\begin{figure}[h]
\centering
\includegraphics[scale=0.22]{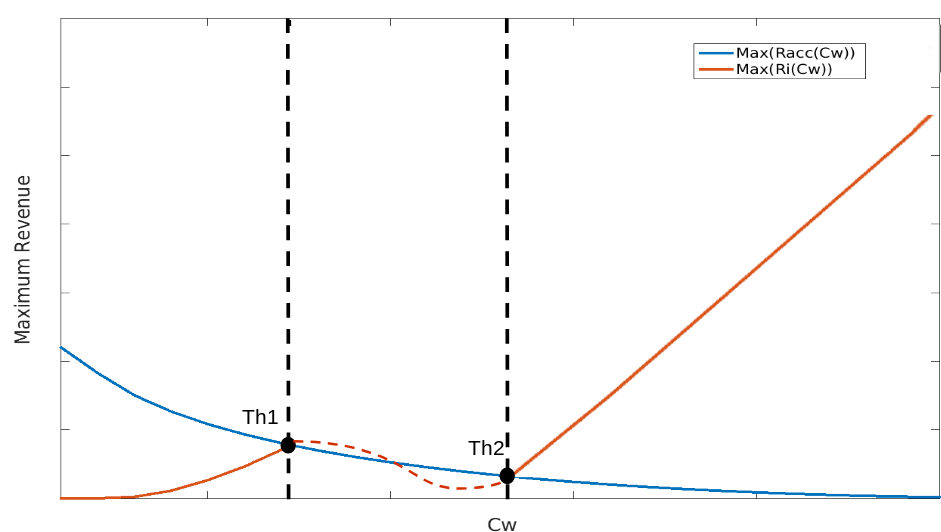} 
\caption{\label{Th:algo} Thresholds Th1 and Th2 as a function of waiting cost}
\end{figure}

Based on the two previous propositions, as a corollary, depending on the waiting cost (the level of sensibility of customers to the waiting time into the queue) the provider has to propose different pricing mechanism in order to maximize his revenue. Specifically, when this sensitivity is low, the provider has an interest to charge the access whereas, when the sensitivity is high, he will prefer to charge the information.\\
Note that there exist at most two thresholds as illustrated on figure \ref{Th:algo} that determine which pricing mechanism is the best. Those thresholds cannot be given in a closed-form but they can be easily determined numerically as illustrated in the numerical section that follows.

%% file: numerical.tex
In this section, we first show the performance of our heuristic when provider charges the queue length information to new customer. Second, when provider charges the access to the queue, we analyze the optimal price $C_{Acc}^*$  depending on the system settings. Finally, we study which pricing mechanism gives the best revenue for the provider depending on the waiting cost. We perform all our simulations using Matlab server R2014b on a Dell Optiplex 9030 i5 running Ubuntu MATE 16.04.
\subsection{Heuristic's performance}
Heuristic approach is scored using different steps (we use this step to update the new revenue value). Table \ref{tabPer} shows the different steps and their total execution time. We provide also the gap between the result obtained using the heuristic and the optimal value given by Matlab solver.

\begin{center}
\begin{tabular}{|*{4}{c|}}
  \hline
 \textbf{Step size}  & 1 & 5.$10^{-1}$ & $10^{-1}$  \\
  \hline
  \textbf{Total time (s)}  & 0.095 & 0.146 & 54.640 \\
  \hline
 \textbf{Gap} & 2.$10^{-4}$ & $10^{-4}$ & $10^{-8}$\\
 \hline
\end{tabular}
\captionof{table}{Performance summary of the heuristic}
\label{tabPer}
\end{center}
The results show clearly that the heuristic plays its role in the optimization of the revenue provider.
The algorithm returns a very good approximate value but with high execution time when we use a small step. With the lower step configured, the result obtained is close to $10^{-8}$ of the exact value. On the other side, when we use higher steps, the gap between the exact value and the value returned by our approach is important but we get the result with a reasonable execution time. 
\vspace{0.5cm}

\subsection{Charging the access to the queue}
\begin{figure}[h]
\centering
\includegraphics[scale=0.14]{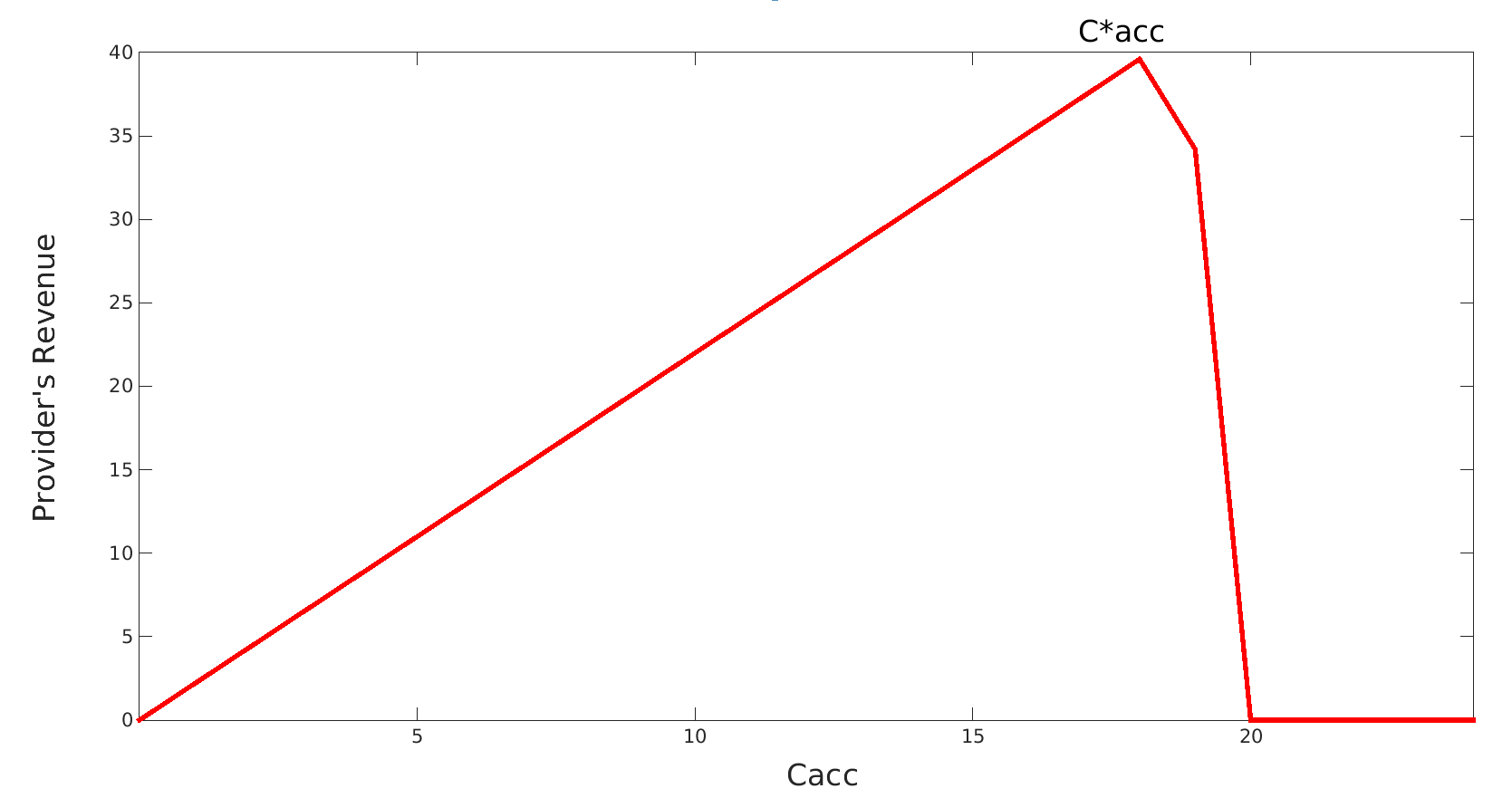} 
\caption{\label{fig:Ra} Revenue as a function of Access fee when: $\mu=2.8$, $\lambda=2.2$, $R=20$ and $C_W=1$}
\end{figure}

\begin{figure}[h]
\centering
\includegraphics[scale=0.14]{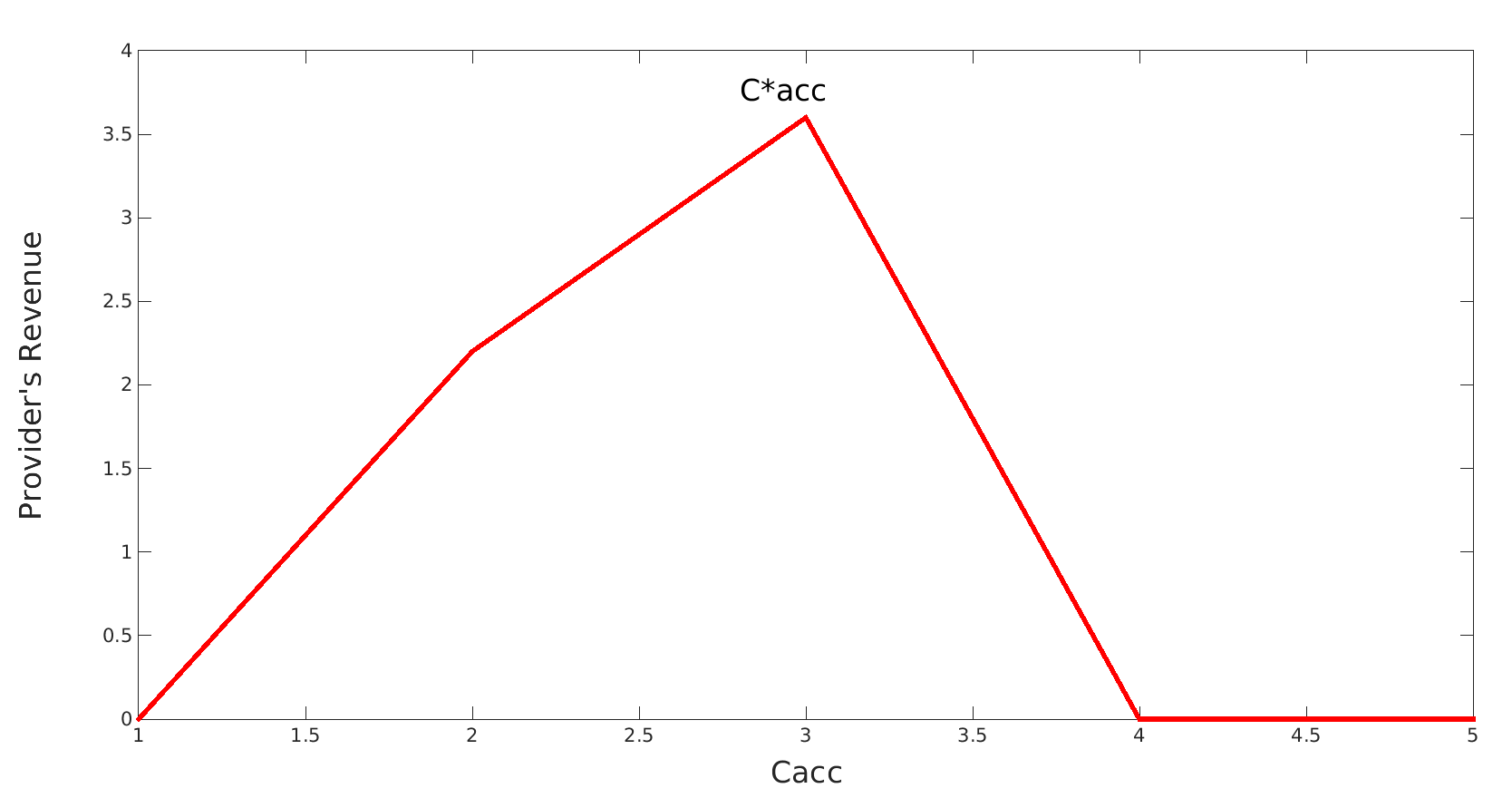} 
\caption{\label{fig:Ra2} Revenue as a function of Access fee when: $\mu=2.8$, $\lambda=2.2$, $R=3$ and $C_W=1$}
\end{figure}

Figure \ref{fig:Ra} and \ref{fig:Ra2} show the provider's revenue as a function of the access fee. Here we use  two system configurations in order to get the optimal price $C_{Acc}^*$. In figure \ref{fig:Ra}, we obtain that  $C_{Acc}^*=(R-\frac{C_W}{\mu-\lambda})=18.33$. This maximum revenue is obtained when all customers join the queue, i.e. $q^*=1$.

Whereas, by considering a lower reward value, we obtain $C_{Acc}^*$ = $(R-\sqrt{R\frac{C_W}{\mu}})=1.96$. This result is illustrated on figure \ref{fig:Ra2} and, in this case, the maximum revenue is obtained when part of customers join the queue, i.e. $0<q^*<1$.

\vspace{0.5cm}
\subsection{Optimal pricing policy}
\begin{figure}[h]
\centering
\includegraphics[scale=0.14]{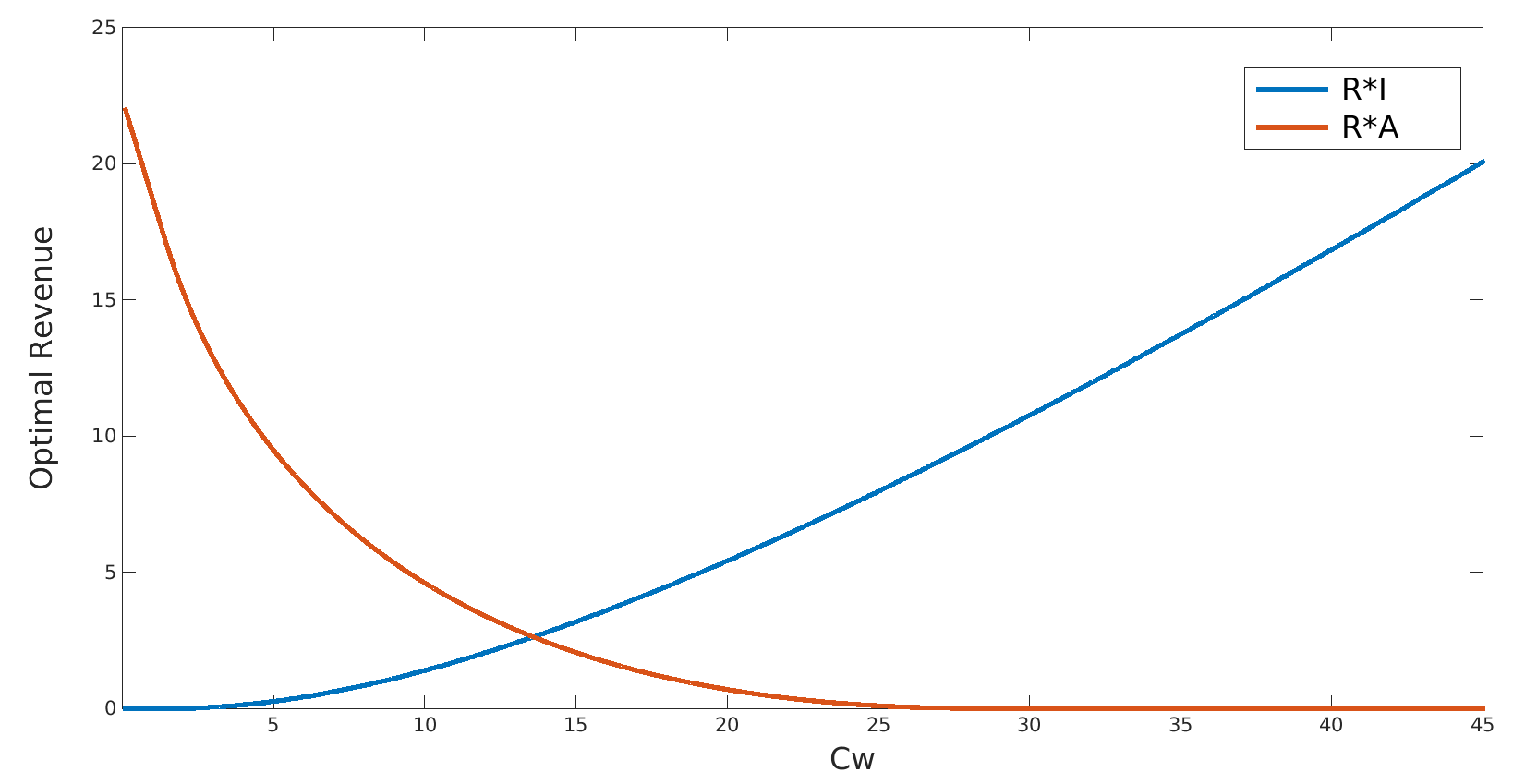} 
\caption{\label{fig:MaxR}Optimal revenues as a function of waiting cost when: $\mu=2.8$, $\lambda=2.2$ and $R=10$}
\end{figure}

Figure \ref{fig:MaxR} gives a numerical example of the maximum revenue obtained by the service provider depending on the pricing mechanism used and as a function of waiting cost $C_W$. We observe that the maximum revenue based on the fee of access is monotonically decreasing: this is because of the increasing cost of the waiting time, which makes customers reluctant to enter; they may be not willing to wait for service avoiding paying too much. The result has been proved theoretically in previous section. We observe moreover that the maximum revenues based on inspection cost is monotonically increasing which means that customers have greater motivation to enter the system as the  cost of waiting increases. Then, this illustrate that there exists in fact a unique threshold that determines which pricing policy is the best. As a result, we obtain a threshold that allows the provider to choose the best pricing policy.



%% file: appendix.tex
\textbf{Proof of observation \ref{ob1}}\\
Let us consider the expected utility for not getting the queue length information as:
\begin{equation*}
\begin{split}
U_{NI}&=\sum_{i=0}^{\infty}\pi_i(R-C_W\frac{ i+1 }{\mu }),\\
&=U_I+C_I+\sum_{i=n_e}^{\infty}\pi_i(R-C_W\frac{i+1}{\mu}).
\end{split}
\end{equation*}
Note that $\forall i \geq n_e$ we have $R<C_W\frac{i+1}{\mu}$ by definition of the threshold $n_e$ given in equation (\ref{Naor}). Then, when $C_I=0$ and for any $p \in [0,1]$, we have:
$$
U_{NI}(p)<U_I(p).
$$
\cqfd
\textbf{Proof of Theorem \ref{theo1}}\\
At equilibrium $p^*$ , we have :
$$
U_{NI} (p^*)= U_I (p^*),
$$
which is equivalent to :
\begin{equation*}
\begin{split}
&\sum_{i=n_e}^{\infty}\pi_i(p^*)(R-C_W\frac{ i+1 }{\mu })=-C_I,\\
&\sum_{i=n_e}^{\infty}\pi_i(p^*)(C_W\frac{ i+1 }{\mu }-R)=C_I,
\end{split}
\end{equation*}
with:
\begin{equation*}
\begin{split}
\pi_i(p)=\frac{\rho^i(1-p)^{i-n_e}}{\frac{1-\rho^{n_e}}{1-\rho}+\frac{\rho^{n_e}}{1-(1-p)\rho}}.
\end{split}
\end{equation*}
We get to solve the following equation :
\begingroup\makeatletter\def\f@size{7.5}\check@mathfonts
$$
\sum_{i=n_e}^{\infty}\rho^i(1-p)^{i-n_e}(C_W\frac{ i+1 }{\mu }-R)=C_I\left(\frac{1-\rho^{n_e}}{1-\rho}+\frac{\rho^{n_e}}{1-(1-p)\rho} \right).
$$
\endgroup
From the left part, we get :
\begingroup\makeatletter\def\f@size{8}\check@mathfonts
\begin{equation*}
\begin{split}
&\rho^{n_e}\sum_{i=0}^{\infty}(\rho(1-p))^{i}(C_W\frac{ i+n_e+1 }{\mu }-R)\\
&=\rho^{n_e}\sum_{i=0}^{\infty}(\rho(1-p))^{i}C_W\frac{ i}{\mu }+(\frac{C_W(n_e+1)}{\mu}-R)\rho^{n_e}\sum_{i=0}^{\infty}(\rho(1-p))^{i}\\
&=\frac{\rho^{n_e}C_W(1-p)\rho}{\mu(1-(1-p)\rho)^2}+\frac{\rho^{n_e}C_W(n_e+1)-R\rho^{n_e}\mu}{\mu(1-(1-p)\rho)}.\\
\end{split}
\end{equation*}
\endgroup
And, from the right part, we get: 
\begingroup\makeatletter\def\f@size{8}\check@mathfonts
\begin{equation*}
\begin{split}
&C_I\left(\frac{1-\rho^{n_e}}{1-\rho}+\frac{\rho^{n_e}}{1-(1-p)\rho} \right)\\
&=\frac{C_I(1-\rho^{n_e})(1-(1-p)\rho)+C_I\rho^{n_e}(1-\rho)}{(1-(1-p)\rho)(1-\rho)}.
\end{split}
\end{equation*}
\endgroup
Then, we obtain:
\begingroup\makeatletter\def\f@size{8}\check@mathfonts
\begin{equation*}
\begin{split}
&\frac{\rho^{n_e}C_W(1-p)\rho+(\rho^{n_e}C_W(n_e+1)-R\rho^{n_e}\mu)(1-(1-p)\rho)}{\mu(1-(1-p)\rho)^2}\\
&=\frac{C_I(1-\rho^{n_e})(1-(1-p)\rho)+C_I\rho^{n_e}(1-\rho)}{(1-(1-p)\rho)(1-\rho)}.
\end{split}
\end{equation*}
\endgroup
Which is equivalent to:
\begingroup\makeatletter\def\f@size{8}\check@mathfonts
\begin{equation*}
\begin{split}
&(1-p)^2 \left( \mu C_I (1-\rho^{n_e})\rho^2\right) -(1-p)( C_I(1-\rho^{n_e})\rho \mu+ K_2\rho \mu\\
&+(1-\rho)(\rho^{n_e+1}C_W- k_1 \rho))+K_2 \mu-K_1(1-\rho)=0.
\end{split}
\end{equation*}
\endgroup
To simplify the following computation, we define : 
$$P=(1-p).$$
So, we get to solve the following second degree equation : 
\begin{equation*}
\begin{split}
&K_3 P^2+K_4 P+ K_5=0.\\
\end{split}
\end{equation*}
If $\Delta >$ 0, the solutions of the second degree equation will be:\\

\begin{center}
$P_1=\left(\frac{-K_4-\sqrt{\Delta}}{2 K_3}\right)$ and
$P_2=\left(\frac{-K_4+\sqrt{\Delta}}{2 K_3}\right).$
\end{center}

Which is equivalent to :
\begin{equation*}
\begin{split}
P_1&=\frac{1}{2\rho}+ \frac{(1-\rho^{n_e+1})}{ (1-\rho^{n_e})2 \rho}+ \frac{(1-\rho)(\rho^{n_e+1}C_W- cst1\rho)}{2\mu C_I (1-\rho^{n_e})\rho^2}\\
&-\frac{\sqrt{\Delta}}{2\mu C_I (1-\rho^{n_e})\rho^2}, \\
P_2&=\frac{1}{2\rho}+ \frac{(1-\rho^{n_e+1})2\rho}{ (1-\rho^{n_e})}+ \frac{(1-\rho)(\rho^{n_e+1}C_W- cst1 \rho)}{2\mu C_I (1-\rho^{n_e})\rho^2}\\
&+\frac{\sqrt{\Delta}}{2 \mu C_I (1-\rho^{n_e})\rho^2}.
\end{split}
\end{equation*}
We wish to determine the sign of $p_2$. We have: \\
\begin{center}
 $\frac{1}{2\rho} > \frac{1}{2}$ because $\rho$ <1 $\Rightarrow$ $\frac{1}{\rho}$ >1,\\
 
 $\frac{(1-\rho^{n_e+1})}{ (1-\rho^{n_e})} \frac{1}{2 \rho}> \frac{1}{2}$ because $(1-\rho^{n_e+1})$ > $(1-\rho^{n_e})$,

 $\frac{\sqrt{\Delta}}{2\mu C_I (1-\rho^{n_e})\rho^2}>0$ because $\sqrt{\Delta}>0$ and $(1-\rho^{n_e})>0$,\\
 
 $\frac{(1-\rho)(\rho^{n_e+1}C_W- K_1 \rho)}{2\mu C_I (1-\rho^{n_e})\rho^2}>0$ because\\
 $\frac{(1-\rho)(\rho^{n_e+1}C_W- K_1\rho)}{2\mu C_I (1-\rho^{n_e})\rho^2}$=$\frac{(1-\rho)\rho^{n_e+1}(R\mu-C_W \lfloor \frac{R\mu}{C_W} \rfloor)}{2 \mu C_I (1-\rho^{n_e})\rho^2} > 0.$
\end{center}
We can deduce that $P_2>$ 1, and as $p = 1-P$, we obtain $p_2< $0. As a result, for each parameters set of the system, a unique equilibrium $p^∗$ can exists defined by: 
$$p^*=p_1=1-P_1$$
or equivalently, 
\begingroup\makeatletter\def\f@size{8}\check@mathfonts
$$
p^*=\frac{2\rho-1}{2\rho}-\frac{(1-\rho^{n_e+1})}{2\rho(1-\rho^{n_e})}-\frac{(1-\rho)(\rho^{n_e+1}C_W-K_1\rho)-\sqrt \Delta}{2\mu C_I(1-\rho^{n_e})\rho^2}.
$$
\endgroup
\cqfd

\textbf{Proof of Proposition \ref{prop2}}\\
We distinguish three cases to obtain the revenue:
\vspace{0.3cm}

$R_A(C_{Acc})$ =\\
\begin{equation*}
\left\{
    \begin{array}{ll}
       \lambda C_{Acc} & \mbox{ if    }  C_{Acc} \le R-\frac{C_W}{\mu-\lambda},  \\
         C_{Acc}( \mu-\frac{C_W}{R-C_{Acc}}) & \mbox{ if    } R-\frac{C_W}{\mu-\lambda}  \le C_{Acc} \le R-\frac{C_W}{\mu},  \\
        0 & \mbox{ if   } C_{Acc} > R-\frac{C_W}{\mu}. \\
    \end{array}
\right.
\end{equation*}
As a result, we obtain optimal access fee $C_{Acc}^*$ in each $C_{Acc}$ interval,  defined as follows: \\
When $ C_{Acc} \le R-\frac{C_W}{\mu-\lambda}$, we have:

        $$C_{Acc}^*=R-\frac{C_W}{\mu-\lambda}.$$
When $R-\frac{C_W}{\mu-\lambda} \le C_{Acc} \le R-\frac{C_W}{\mu}$, we have:
\begin{equation*}
\begin{split}
R_A(C_{Acc})&=\lambda q^*(C_{Acc}) C_{Acc},\\
&= C_{Acc}\left( \mu-\frac{C_W}{R-C_{Acc}} \right),\\
&=C_{Acc}\mu-\frac{C_{Acc}C_W}{R-C_{Acc}}.
\end{split}
\end{equation*}
This function is $C^{\infty}$  and taking its derivative equals to zero gives:\\
\begin{equation*}
\begin{split}
R_A'(C_{Acc})=0 &\Leftrightarrow \mu-\frac{C_W(R-C_{ACC})+C_{ACC} C_W}{(R-C_{ACC})^2}=0\\
&\Leftrightarrow \mu-\frac{R C_W}{(R-C_{ACC})^2}=0\\
&\Leftrightarrow (R-C_{Acc})^2=\frac{R C_W}{\mu}\\
&\Leftrightarrow C_{Acc}^*=R-\sqrt{\frac{R C_W}{\mu}}. \\
\end{split}
\end{equation*}
When $C_{Acc} > R-\frac{C_W}{\mu} $, we have: $C_{Acc}^*=0$.
\vspace{1cm}\\
Then, the optimal access fee $C_{Acc}^*$ is either: 
\begin{center}
 $R-\frac{C_W}{\mu-\lambda}$ or $R-\sqrt{R\frac{C_W}{\mu}} .$
\end{center}


\cqfd


\textbf{Proof of Proposition \ref{prop3}}\\
When $ C_{Acc} \le R-\frac{C_W}{\mu-\lambda}$, we have the following equation of revenue as a function of waiting cost: \\
\begin{equation*}
\label{Ra1}
\begin{split}
R_A(C_{Acc}(C_W),C_W)&=\lambda q^*(C_{Acc}^*(C_W))C_{Acc}^*(C_W),\\
&=\left(R-\frac{C_W}{\mu-\lambda} \right)\lambda,\\
&=R\lambda- \frac{\lambda C_W}{\mu-\lambda}:=\check{R}(C_W).\\
\end{split}
\end{equation*}
The derivate of this equation gives: \\
\begin{equation*}
\check{R}'(C_W)=\frac{-\lambda}{\mu-\lambda}<0.\\
\end{equation*}
When $R-\frac{C_W}{\mu-\lambda} \le C_{Acc} \le R-\frac{C_W}{\mu}$, we have the following equation of revenue as a function of waiting cost:\\
\begin{equation*}
\label{Ra2}
\begin{split}
R_A(C_{Acc}(C_W),C_W)&=\lambda q^*(C_{Acc}^*(C_W))C_{Acc}^*(C_W), \\
&=\left(\mu-\frac{C_W}{\sqrt{\frac{R C_W}{\mu}}} \right) \left(R-\sqrt{\frac{R C_W}{\mu}} \right),\\
&=C_W-\sqrt{C_W}(R\sqrt{\frac{\mu}{R}}-\mu\sqrt{\frac{R}{\mu}})+\mu R,\\
&=C_W-2\sqrt{\mu R}\sqrt{C_W}+\mu R:=\check{R}(C_W).
\end{split}
\end{equation*}
The derivate of this equation gives: \\
\begin{equation*}
\check{R}'(C_W)= 1-\frac{\sqrt{\mu R}}{\sqrt{C_W}}.\\
\end{equation*}
\begin{equation*}
\begin{split}
\check{R}'(C_W)<0 &\Leftrightarrow 1<\frac{\sqrt{\mu R}}{\sqrt{C_W}} \Leftrightarrow \sqrt{C_W}<\sqrt{\mu R}
\end{split}
\end{equation*}
When, 
\begin{equation*}
\begin{split}
&C_{Acc}> 0 \Leftrightarrow R>\sqrt{\frac{R C_W}{\mu}}\Leftrightarrow \sqrt{R \mu} > \sqrt{C_W}\\
\end{split}
\end{equation*}
\cqfd


\textbf{Proof of Proposition \ref{prop4}}\\
We have the following equation of revenue as function of the waiting cost:
\begin{equation*}
    \begin{split}
        R_I(C_I(C_W), C_W) =\lambda  p^*(C_I^*(C_W)) C_I^*(C_W).
    \end{split}
\end{equation*}
When $C_W$ tends to zero and for any $C_I$, we have: \\
$$U_I=\sum_{i=0}^{n_e-1}  \pi_i(R-C_W\frac{ i+1 }{\mu })- C_I \underset{C_W \to 0}{\longrightarrow} R \sum_{i=0}^{n_e-1}(\pi_i)- C_I,$$
$$U_{NI}=\sum_{i=0}^{\infty}\pi_i(R-C_W\frac{ i+1 }{\mu })\underset{C_W \to 0}{\longrightarrow} R.$$
As a result, we obtain: \\
$$R_I(C_I(C_W), C_W) \underset{C_W \to 0}{\longrightarrow}0.$$
Because,
$$U_{NI}>U_I \Longrightarrow p^*=0. $$
When $C_W$ tends to infinity and for any $C_I$, we have: \\
$$U_{NI}=U_I+\sum_{i=n_e}^{+\infty1}  \pi_i(R-C_W\frac{ i+1 }{\mu })+ C_I \underset{C_W \to 0}{\longrightarrow} -\infty, $$ 
$$U_{NI}-U_I< 0 \Longrightarrow p^*=1.$$ 
As a result, we obtain:\\
$$ R_I(C_I(C_W), C_W) \underset{C_W \to +\infty}{\longrightarrow} +\infty.$$
\cqfd